\documentclass[a4paper,10pt]{amsart}

\newtheorem{theo}{Theorem}[section]

\newtheorem{prop}[theo]{Proposition}
\newtheorem{lem}[theo]{Lemma}
\newtheorem{cor}[theo]{Corollary}

\newtheorem{rema}[theo]{Remark}
\newtheorem{remas}[theo]{Remarks}

\def \dem {\paragraph{ \em Proof. }}

\def \Romannumeral #1 {\expandafter\uppercase\expandafter {\romannumeral #1} }
\def \br {{\rm{Br \,}}}

\def \P {{\bf P}}
\def \tors {{\rm tors}}

\def \nr {{\rm nr}}

\def \pic {{\rm {Pic\,}}}
\def \Div {{\rm{Div\,}}}
\def \gal {{\rm{Gal\,}}}

\def \calo {{\mathcal O}}
\def \T {{\mathcal T}}

\def \calf {{\mathcal F}}

\def \calt {{\mathcal T}}
\def \spec {{\rm{Spec\,}}}

\def \Hom {{\rm {Hom}}}

\def \nr {{\rm nr}}

\def \Z {{\bf Z}}
\def \Q {{\bf Q}}

\def \RR {{\bf R}}

\def \im {{\rm {Im}}}
\def \ker {{\rm {Ker}}}
\def \G {{\bf G}_m}
\def \A {{\bf A}}

\def\smallsquare{\vbox{\hrule\hbox{\vrule height 1 ex\kern 1 ex\vrule}\hrule}}
\def\enddem{\hfill \smallsquare\vskip 3mm}

\def \abstract{\paragraph{Abstract. }}
\usepackage{amscd}
\usepackage{amssymb}
\usepackage{amsmath}

\everymath{\displaystyle}

\title{Weak approximation for tori over $p$-adic function fields}
\author{David Harari, Claus Scheiderer and Tam\'as Szamuely}
\address{Universit\'e de Paris-Sud Math\'ematique, B\^atiment 425, 91405 Orsay, France}
\email{David.Harari@math.u-psud.fr}
\address{Fachbereich Mathematik and Statistik, Universit\"at Konstanz, D--78457 Konstanz, Germany}
\email{claus.scheiderer@uni-konstanz.de}
\address{Alfr\'ed R\'enyi Institute of Mathematics, Hungarian Academy of Sciences, Re\'altanoda utca 13--15, H-1053 Budapest, Hungary and Central European University, N\'ador utca 9, H-1051 Budapest, Hungary}
\email{szamuely.tamas@renyi.mta.hu}

\date{\today}

\def \nr {{\rm nr}}

\def \tors {{\rm tors}}

\DeclareFontFamily{U}{wncy}{}
\DeclareFontShape{U}{wncy}{m}{n}{%
   <5>wncyr5%
   <6>wncyr6%
   <7>wncyr7%
   <8>wncyr8%
   <9>wncyr9%
   <10>wncyr10%
   <11>wncyr10%
   <12>wncyr6%
   <14>wncyr7%
   <17>wncyr8%
   <20>wncyr10%
   <25>wncyr10}{}
\DeclareMathAlphabet{\cyrille}{U}{wncy}{m}{n}
\def\Sha{\cyrille X}
\def \R{{\bf R}}

\usepackage{graphicx}

\makeatletter
\DeclareRobustCommand\widecheck[1]{{\mathpalette\@widecheck{#1}}}
\def\@widecheck#1#2{%
    \setbox\z@\hbox{\m@th$#1#2$}%
    \setbox\tw@\hbox{\m@th$#1%
       \widehat{%
          \vrule\@width\z@\@height\ht\z@
          \vrule\@height\z@\@width\wd\z@}$}%
    \dp\tw@-\ht\z@
    \@tempdima\ht\z@ \advance\@tempdima2\ht\tw@ \divide\@tempdima\thr@@
    \setbox\tw@\hbox{%
       \raise\@tempdima\hbox{\scalebox{1}[-1]{\lower\@tempdima\box
\tw@}}}%
    {\ooalign{\box\tw@ \cr \box\z@}}}
\makeatother

\begin{document}
\maketitle

\setcounter{section}{-1}

\noindent{\small{\bf Abstract.} We study local-global questions for Galois cohomology over the function field of a curve defined over a p-adic field, the main focus being weak approximation of rational points. We construct a 9-term Poitou--Tate type exact sequence for tori over a field as above (and also a 12-term sequence for finite modules). Like in the number field case, part of the sequence can then be used to analyze the defect of weak approximation for a torus. We also show that the defect of weak approximation is controlled by a certain subgroup of the third unramified cohomology group of the torus.}\bigskip

\section{Introduction}

This paper is the companion piece to \cite{pchp}, containing investigations concerning local-global questions for tori defined over the function field $K$ of a curve over a finite extension of $\Q_p$. As recalled in the introduction of \cite{pchp}, our project has been motivated by the recent awakening of
 interest in local-global principles for group schemes defined over fields of cohomological dimension strictly greater than 2, as documented in the work of Harbater, Hartmann and Krashen \cite{hhk} as well as Colliot-Th\'el\`ene, Parimala and Suresh (\cite{ctps}, \cite{ctps2}).

In \cite{pchp} two of us are studying the Hasse principle for torsors under tori over a field $K$ as above, the main tool being a global duality theorem. Presently our main concern is weak approximation for tori. As in the classical case over number fields, this necessitates going beyond duality and establishing a Poitou--Tate type exact sequence for tori over $K$. Here is our first main result:

\begin{theo}\label{varpoitores} {\rm (= Theorem \ref{poitores})}
Let $k$ be a finite extension of $\Q_p$, and
$T$ a torus defined over the function field of a smooth proper $k$-curve $X$. There is an exact
sequence of topological groups
$$
\begin{CD}
0 @>>> H^0(K,T)_{\wedge}  @>>> \P^0(T)_{\wedge}  @>>> H^2(K,T')^D \cr
&&&&&& @VVV \cr
&& H^1(K,T') ^D @<<< \P^1(T) @<<< H^1(K,T) \cr
&& @VVV \cr
&& H^2(K,T) @>>> \P^2(T)_{\tors} @>>> (H^0(K,T')_{\wedge})^D @>>> 0.
\end{CD}
$$
\end{theo}

In this sequence, $T'$ denotes the dual torus of $T$ (i.e. the torus whose character group is the cocharacter group of $T$) and the groups $\P^i(T)$ are certain restricted products of local cohomology groups to be defined in Section \ref{pt}. Furthermore, the subscript `tors' stands for the torsion subgroup,  the subscript $\wedge$ indicates the inverse limit of mod $n$ quotients for all $n>0$, and the superscript $D$ means continuous dual with values in $\Q/\Z$. In Section \ref{pt} we shall also construct a 12-term exact sequence of similar type for finite group schemes over $K$.

As in the number field case, part of the sequence can be used to analyze the defect of weak approximation for a torus. The question here is whether for $T$ as above the group of points $T(K)$ is dense in the topological product of the groups $T(K_v)$, where $K_v$ denotes the completion of $K$ with respect to the discrete valuation associated with a closed point $v\in X$, equipped with its natural topology. The answer is yes when the torus is $K$-rational, but in general there is an obstruction, as shown by the following analogue of a classical result of Voskresenskii (\cite{vosk}; see also \cite{requ} and \cite{sansuc}).

\begin{theo}\label{varmainwa} {\rm (= Theorem \ref{mainwa} $b)$)} For $K$ and $X$ as above, denote by $X^{(1)}$ the set of closed points of $X$. For a $K$-torus $T$ let $\overline{T(K)}$ denote the closure of
$T(K)$ in the topological product of the $T(K_v)$ for all $v\in X^{(1)}$. There is an exact sequence

$$
0 \to \overline{T(K)} \to \prod_{v \in X^{(1)}} T(K_v) \to
\Sha^2_{\omega}(T')^D \to \Sha^1(T) \to 0,
$$
where $\Sha^1(T)\subset H^1(K,T)$ is the subgroup of classes whose image in $H^1(K_v, T)$ is trivial for all $v\in X^{(1)}$, and $\Sha^2_{\omega}(T')\subset H^2(K, T')$ is the subgroup of classes whose image in $H^2(K_v, T')$ is trivial for all but finitely many $v$.
\end{theo}

Note that, in contrast to the number field case, in our situation the group $\Sha^2_\omega(T')$ may be infinite (see Proposition \ref{exsha} below), whereas $\Sha^1(T)$ is always finite.

In the last section we reinterpret the above theorem using the reciprocity obstruction to weak approximation first studied by J-L. Colliot-Th\'el\`ene. In our case this obstruction is defined by means of the subgroup $H^3_{\rm nr}(K(T), \Q/\Z(2))\subset H^3(K(T), \Q/\Z(2))$  of elements having trivial residues for all discrete valuations of $K(T)$ trivial on $K$. (Here $\Q/\Z(2)$ denotes the direct limit of the torsion \'etale sheaves $\mu_n^{\otimes 2}$ for all $n$.) As we shall recall in detail in Section \ref{recwa}, evaluating classes at $K_v$-points of $T$ induces a pairing
$$
\prod_{v\in X^{(1)}}T(K_v)\times H^3_{\rm nr}(K(T), \Q/\Z(2))\to\Q/\Z
$$
which annihilates the closure of the diagonal image of $T(K)$  on the left. We shall prove:

\begin{theo}\label{varmainrecwa} {\rm (= Theorem \ref{mainrecwa})}
For a $K$-torus $T$ each system $(P_v)$ of local points of $T$ lying in the left kernel in the above pairing is in the closure of the diagonal image of $T(K)$ for the product topology.
\end{theo}

In fact, in Theorem \ref{mainrecwa} we shall prove slightly  more: we construct a homomorphism
$u: \Sha^2_{\omega}(T') \to H^3_{\rm nr}(K(T),\Q/\Z(2))$ relating the obstruction of Theorem \ref{varmainwa} to the group $H^3_{\rm nr}(K(T),\Q/\Z(2))$, and
show that it is enough to consider the restriction of the pairing to $\im (u)$ in order to obtain the conclusion of the theorem. \smallskip

The above theorems were announced, in slightly weaker form, by the second
author in seminar lectures given in 2002. The proofs given in the present
text were worked out by the first and third authors in 2012-13, in
connection with the results of \cite{pchp}. We are very grateful to Jean-Louis Colliot-Th\'el\`ene and Yves de Cornulier for instructive discussions. Part of this work has been done while the first and third  authors were visiting each other's home institutions whose hospitality was greatly appreciated. The first author would also like to thank the Centre Bernoulli (Lausanne) for hospitality and the third author the Hungarian Science Foundation OTKA for partial support under grant No. NK81203.

\bigskip

\noindent{\bf Notation and conventions.} The following conventions will be adopted throughout the article.\smallskip

\noindent{\em Abelian groups.\/}  Given an abelian group $A$, we shall denote by $\overline A$
the quotient of $A$ by its maximal divisible subgroup $\Div A$. The $n$-torsion subgroup and (for $\ell$ prime) the $\ell$-primary torsion subgroup of  $A$ will be denoted by $_n A$ and $A \{ \ell\}$, respectively.
If $A$ is torsion, it is said to be {\it of cofinite type} if for each
$n >0$, the subgroup $_n A$ is finite. For $A$ a topological abelian group the notation $A^D$ will stand for the group of continuous homomorphisms $A\to\Q/\Z$.  The functor $A \mapsto A^D$
is an anti-equivalence of categories between $\ell$-primary torsion
abelian groups of cofinite type (with the discrete topology) and $\Z_{\ell}$-modules of finite type (with the profinite topology).

The notation $A^\wedge$ will stand for the profinite completion of $A$, i.e. the inverse limit of its finite quotients. We shall denote by $A_\wedge$ the inverse limit of the quotients $A/nA$ for all $n>0$.

Unless otherwise stated, all cohomology groups will be taken with respect to the \'etale topology.\smallskip

\noindent {\em Tori.\/} For a torus $T$ over a field $F$ we shall denote by $\widehat T$ its character module and by $\widecheck T$ its module of cocharacters. These are finitely generated free $\Z$-modules equipped with a Galois action, and moreover $\widecheck T$ is the $\Z$-linear dual of $\widehat T$. The {\em dual torus} $T'$ of $T$ is by definition the $F$-torus with character group $\widecheck T$. The torus $T$ is {\em quasi-trivial} if $\widehat T$ is a permutation module, i.e. it has a Galois invariant basis. The torus  $T$ is  {\em flasque} if it is split by a finite \'etale Galois cover of the base with group $G$ and $H^{-1}(H, \widehat T)=0$ for all subgroups $H\subset G$. It is {\em coflasque} if it has a similar property for $H^1$ instead of $H^{-1}$. The dual of a flasque torus is coflasque and vice versa. Every torus $T$ has a {\em flasque resolution}, i.e. an exact sequence
\begin{equation}\label{flasque}
1 \to S \to P \to T \to 1
\end{equation}
with $P$ a quasi-trivial and $S$ a flasque torus. This holds in fact not just over fields but over arbitrary bases, under very mild assumptions (\cite{ctflasq}, Theorem 1.3).\smallskip

\noindent{\em Motivic complexes.\/} For $i\geq 0$ and a separated scheme $V$ of finite type over a field $F$ we denote by $z^i(X,\bullet)$ Bloch's cycle complex defined in \cite{hcg}. If $V$ is smooth, the \'etale motivic complex $\Z(i)$ over $V$ is defined as the complex of sheaves $z^i(-,\bullet)[-2i]$ on the small \'etale site of $V$. For $i=1$ we have a quasi-isomorphism of complexes of sheaves $\Z(1)\cong \G[-1]$. The analogously defined Zariski motivic complex will be denoted by $\Z(i)_{\rm Zar}$. For an abelian group $A$ we denote by $A(i)$ the complex $A\otimes\Z(i)$ and similarly for $\Z(i)_{\rm Zar}$ (since the terms of $\Z(i)$ are torsion free, this is the same as the derived tensor product). For $m$ invertible in $F$ we have a quasi-isomorphism of complexes of \'etale sheaves
\begin{equation}\label{Z(i)modm}
\Z/m\Z(i)\stackrel\sim\to\mu_m^{\otimes i}
\end{equation} where $\mu_m$ is the \'etale sheaf of $m$-th roots of unity placed in degree 0 (Geisser--Levine \cite{gl}, Theorem 1.5). Thus we shall also use the notation $\Q/\Z(i)$ for the direct limit of the sheaves $\mu_m^{\otimes i}$ for all $m>0$. \smallskip

\noindent{\em Function fields.\/} Throughout the paper, $k$ will be a finite extension of $\Q_p$ for a prime number $p$, and $X$ a smooth proper geometrically integral curve over $k$. The set of all closed points of $X$ will be denoted by $X^{(1)}$, and its function field by $K$. For a closed point $v\in X^{(1)}$ we denote by $\kappa(v)$ its residue field, and by $K_v$ the completion of $K$ for the discrete valuation induced by $v$. Then $\kappa(v)$ is also the residue field of the ring of integers $\calo_v$ and is a finite extension of $k$. Therefore $K_v$ is a {\em 2-dimensional local field\/}, i.e. a field complete with respect to a discrete valuation whose residue field is a classical local field.

\section{Duality results: summary and complements}\label{one}

In this section we recall the main duality results established in \cite{pchp} and prove some complements that will be needed for the construction of the Poitou--Tate exact sequence.

Let us start with the local theory. Given a discrete valuation $v$ of the function field $K$ coming from a closed point of the curve $X$, the completion of $K$ at $v$ is a 2-dimensional local field. For such a $K_v$ there is a canonical isomorphism $$H^4(K_v, \Z(2))\cong\Q/\Z$$
by \cite{pchp}, Lemma 2.1.

If $T$ is a torus over $K$ with dual torus $T'$, there is a canonical pairing
$$
T\otimes^{\bf L} T'\to\Z(2)[2]
$$
in the derived category of \'etale sheaves over $K_v$. It can be constructed using the quasi-isomorphisms
\begin{equation}\label{dualtordef} T\cong \widecheck T\otimes^{\bf L} \Z(1)[1], \quad
T'= \widehat T\otimes^{\bf L} \Z(1)[1]
\end{equation}
(which follow from the quasi-isomorphism $\Z(1)[1]\cong\G$) by putting the pairings $\widecheck T\otimes\widehat T\to\Z$ and $\Z(1)[1]\otimes^{\bf L}\Z(1)[1]\to\Z(2)[2]$ together.

We therefore obtain cup-product pairings

\begin{equation}\label{localpairing}
H^i(K_v, T)\otimes H^{2-i}(K_v, T')\to H^4(K_v, \Z(2))\cong\Q/\Z
\end{equation}
for $i=0,1,2$. Concerning these pairings, we have:

\begin{prop} \label{localdual}
The pairing $(\ref{localpairing})$ is a perfect pairing of finite groups for $i=1$. For $i=0$ it becomes a perfect pairing between a profinite and a torsion group, after replacing $H^0(K_v, T)$ by its profinite completion.
\end{prop}

\begin{dem} See \cite{pchp}, Proposition 2.2.
\end{dem}

Recall also (e.g. from \cite{adt}, I. 2.17) that for a finite $K_v$-group scheme $F$
there are perfect pairings of finite groups
\begin{equation}\label{finitepairing}
H^i(K_v, F)\otimes H^{3-i}(K_v, F')\to H^3(K_v, \Q/\Z(2))\cong\Q/\Z
\end{equation}
where $F':=\Hom(F,\Q/\Z(2))$. Note that in the case when $F$ is the $m$-torsion subgroup ${}_mT$ of a $K$-torus $T$ for some $m>0$, we have $F'\cong {}_mT'$ in view of the chain of isomorphisms
$$
\Hom({}_mT,\mu_m^{\otimes 2})\cong \Hom(\widecheck T\otimes\mu_m,\mu_m^{\otimes 2})\cong \Hom(\widecheck T,\mu_m)\cong {}_mT'.
$$

We now establish some additional properties of local duality for finite group schemes and for tori.

\begin{prop} \label{annihil}
Assume
$F$ is a finite group scheme over $K_v$ that extends to a finite and \'etale group scheme $\calf$ over
$\spec \calo_v$. Then under the pairing
$$H^1(K_v,F) \times H^2(K_v,F') \to \Q/\Z  $$
the subgroups $H^1(\calo_v,\calf) \subset H^1(K_v,F)$ and
$H^2(\calo_v,\calf') \subset H^2(K_v,F')$ are exact annihilators of each other.
\end{prop}

In fact, it is not obvious a priori that the maps $H^1(\calo_v,\calf) \to H^1(K_v,F)$ and
$H^2(\calo_v,\calf') \to H^2(K_v,F')$ are injective, but this will follow from exact sequences (\ref{eqnalihil1}) and (\ref{eqnalihil2}) in the proof below.\medskip

\dem First of all, the cup-product pairing
\begin{equation}\label{cupov}
H^i(\calo_v,\calf) \times H^{3-i}(\calo_v,\calf')\to H^3(\calo_v,\Q/\Z(2))
\end{equation}
induced by $\calf\otimes \calf'\to\Q/\Z(2)$ is trivial for all $i$, since $$H^3(\calo_v,\Q/\Z(2))\cong H^3(\kappa(v),\Q/\Z(2))=0,$$ the residue field $\kappa(v)$ being of cohomological dimension 2.

Let now $I_v$ be the inertia subgroup of the Galois group $G_v$ of $K_v$, and denote by $G(v)=G_v/I_v$ the
Galois group of the residue field $\kappa(v)$. By assumption, we may identify $F$ with a $G_v$-module with trivial $I_v$-action, and we have
$H^i(\calo_v,\calf)=H^i(G(v),F)$ (and likewise for $F'$).
Since $F$ is unramified, the Hochschild--Serre spectral sequence yields a short exact sequence
\begin{equation}\label{eqnalihil1}
0 \to H^1(G(v),F) \to H^1(G_v,F) \to H^0(G(v),H^1(I,F))
\to 0.
\end{equation}
Using moreover the fact that the exact sequence of Galois
groups
$$1 \to I_v \to G_v \to G(v) \to 1$$
splits (\cite{cogal}, \S II.4.3), we similarly get a short exact sequence
\begin{equation}\label{eqnalihil2}
0 \to H^2(G(v),F') \to H^2(G_v,F') \to H^1(G(v),H^1(I,F'))
\to 0.
\end{equation}
On the other hand, the $G(v)$-module  $H^1(I,F)=\Hom(I,F)$ is isomorphic to
$F(-1)$, which is the module of characters of $F'$
(and vice versa).
Therefore by Tate's local duality theorem for finite Galois modules over $p$-adic fields, the last terms in exact sequences (\ref{eqnalihil1}) and (\ref{eqnalihil2}) are dual finite abelian groups. It follows that $H^1(K_v,F)/H^1(\calo_v,\calf)$ and
$H^2(K_v,F')/H^2(\calo_v,\calf')$ are finite of the same cardinality, which shows that
$H^1(\calo_v,\calf)$ and $H^2(\calo_v,\calf')$ are indeed exact annihilators.

\enddem

\begin{prop}\label{annihiltori}
Let  $T$
be a torus over $K_v$ that extends to a torus $\calt$ over $\spec \calo_v$.
Then~:

\smallskip

a) Under the pairing
$$H^1(K_v,T) \times H^1(K_v,T') \to \Q/\Z  $$
the subgroups $H^1(\calo_v,\calt) \subset H^1(K_v,T)$ and
$H^2(\calo_v,\calt') \subset H^2(K_v,T')$ are exact annihilators of each other.

\smallskip

b) The annihilator of $H^2(\calo_v,\calt')$ under the pairing
$$H^0(K_v,T)^{\wedge} \times H^2(K_v,T') \to \Q/\Z  $$
is $H^0(\calo_v,\calt)^{\wedge}$.
\end{prop}

\noindent Note that for each $n >0$ the groups
$H^0(K_v,T)/n$ and their subgroups $H^0(\calo_v, \calt)/n$
are finite. Hence the profinite completions $H^0(K_v,T)^{\wedge}$ and
$H^0(\calo_v, \calt)^{\wedge}$ equal the `$n$-adic completions'
$H^0(K_v,T)_{\wedge}$ and $H^0(\calo_v, \calt)_{\wedge}$, respectively.

Also, as in the previous proposition, the injectivity of the maps $H^1(\calo_v,\calt) \to H^1(K_v,T)$ and
$H^2(\calo_v,\calt') \to H^2(K_v,T')$ will be justified in the course of the proof.\medskip

\begin{dem} $a)$ As before, the pairing $H^1(\calo_v,\calt)\times H^1(\calo_v, \calt')\to\Q/\Z$ is trivial because $H^3(\calo_v,\Q/\Z(2))=0$,
so it remains to prove that each $t \in H^1(K_v,T)$ orthogonal to $H^1(\calo_v, \calt')$ comes from $H^1(\calo_v,\calt)$.
By functoriality, for each $n>0$ the image of $t$ in
$H^2(K_v,\, _n T)$ is orthogonal to $H^1(\calo_v,\, _n T')$, hence comes
from $H^2(\calo_v ,\, _n \calt)$ by Proposition~\ref{annihil}.
There is a commutative diagram with exact rows
$$
\begin{CD}
&& H^1(\calo_v,\calt) @>>> H^2(\calo_v,\,{}_n\calt) @>>> H^2(\calo_v,\calt) \\
&& @VVV @VVV @VVV \\
 H^1(K_v,T) @>n>> H^1(K_v,T) @>>>H^2(K_v,\,{}_nT) @>>> H^2(K_v,T).
\end{CD}
$$
Now observe that the map
$H^2(\calo_v,\calt) \to H^2(K_v,T)$ is injective. Indeed, we may identify it with a map $H^2(G/I,T(\overline K_v)^I)\to H^2(G,T(\overline K_v))$ in a Hochschild--Serre spectral
sequence, where $I$ is the inertia subgroup in $G=\gal(\overline K_v|K_v)$;  this map is injective because $H^0(G/I, H^1(I,T(\overline K)))=0$ by Hilbert's Theorem 90 as $T$ is split by an unramified extension.

A diagram
chase now shows that there exists $t_0 \in H^1(\calo_v,\calt)$ such
that $t-t_0 \in nH^1(K_v,T)$. Since this holds for every $n >0$,
 the image of $t$ in the finite group
$H^1(K_v,T)/H^1(\calo_v,\calt)$ is divisible, hence  zero.

\smallskip

$b)$ We already know that the pairing is trivial on
$H^0(\calo_v,\calt)^{\wedge} \times H^2(\calo_v,\calt')$. Fix $n>0$ and assume
$t \in H^0(K_v,T)/n$  is orthogonal to
$_n H^2(\calo_v,\calt')$. By functoriality the image $s$ of $t$ in
$H^1(K_v, _n T)$ is orthogonal to $H^2(\calo_v,_n \calt')$, hence
$s \in H^1(\calo_v, _n \calt )$ by Proposition~\ref{annihil}.
The exact commutative diagram
$$
\begin{CD}
0 @>>> H^0(\calo_v,\calt)/n @>>> H^1(\calo_v,\,{}_n\calt) @>>> H^1(\calo_v,\calt) \\
&& @VVV @VVV @VVV \\
0 @>>> H^0(K_v,T)/n  @>>> H^1(K_v,_n T) @>>>H^1(K_v,T).
\end{CD}
$$
then implies $t \in H^0(\calo_v,T)/n$ because the right vertical map is injective, again by a Hochschild--Serre argument. This proves the statement.
\end{dem}

We conclude this section by recalling the global duality theorems of \cite{pchp}. Given a torus $T$ over our function field $K$, we set
$$\Sha^i(T):=\ker (H^i(K,T) \to \prod_{v \in X^{(1)}} H^i(K_v,T))$$
for $i\geq 1$. Similar notation will be used for a finite group scheme $F$ over $K$.

\begin{theo}\label{dualsha} Given a $K$-torus $T$, there is a perfect pairing of finite
groups
$$\Sha^1(T) \times \Sha^2(T') \to \Q/\Z.$$

Similarly, for a finite group scheme $F$ over $K$ there are perfect pairings of finite
groups
$$\Sha^i(F) \times \Sha^{4-i}(F') \to \Q/\Z$$
for $i=1,2$.
\end{theo}

\begin{dem}
See \cite{pchp}, Theorems 4.1 and 4.4.
\end{dem}

\section{Poitou-Tate exact sequences}\label{pt}

In this section we construct Poitou--Tate type sequences for finite modules and
tori over $K$. We start with the case of finite modules.
Let $F$ be a finite group scheme over $K$. We know that $F$ extends to
a finite \'etale group scheme $\calf$ over a suitable Zariski open $U_0\subset X$. For each $i\geq 0$
denote by
$\P^i(F)$ the topological restricted  product of the $H^i(K_v,F)$  for all $v \in X^{(1)}$
with respect to the images of the maps $H^i(\calo_v,\calf)\to H^i(K_v,F)$ for $v\in U_0$ (we shall use the same notation for a $K$-torus
instead of a finite module). One sees that $\P^i(F)$ does not depend on the choices of $U_0$ and of $\calf$.
We have $$\P^0(F)=\prod_{v \in X^{(1)}} H^0(K_v,F).$$ Also,
$$\P^3(F)=\bigoplus_{v \in X^{(1)}} H^3(K_v,F)$$ since
$$H^3(\calo_v,\calf)=H^3(\kappa(v),F_{\kappa(v)})=0$$
for $v \in U_0$ because $\kappa(v)$ is a $p$-adic field (hence
of cohomological dimension $2$). The group $\P^0(F)$ is profinite (hence compact) being the product of finite groups, and $\P^3(F)$ is discrete as a direct sum of such, but the other two groups are only locally compact.

\smallskip

The following result gives pieces of the Poitou--Tate exact sequence.

\begin{prop} \label{partpoitou}
For $i=1,2,3$ there are exact sequences
\begin{equation} \label{poitou1}
 H^i(K,F) \to \P^i(F) \stackrel{\theta}{\to} H^{3-i}(K,F')^D
\end{equation}
where the map $\theta$ is defined via the local pairings $ \langle \,\, ,\,\, \rangle _v$
by the formula
$$\theta((f_v))(f')=\sum_{v \in X^{(1)}} \langle f_v,f'_v \rangle _v$$
for $(f_v) \in \P^i(F)$ and $f' \in H^{3-i}(K,F')$.
\end{prop}

Note that since for almost all $v$ the module $F$ is unramified at $v$ and moreover
$f_v \in H^i(\calo_v,\calf)$ and $f'_v \in H^{3-i}(\calo_v,\calf')$, the sum
in the above formula is actually finite because the pairings (\ref{cupov}) over $\calo_v$ are trivial.
\medskip

For the proof we need the following lemma (which is in fact valid over an arbitrary base field).

\begin{lem} \label{purlemma}
 Let $V \subset U$ be nonempty open subsets of $U_0$, and let
$\alpha \in H^i(V,\calf)$ be such that the localisation $\alpha_v
\in H^i(K_v,F)$ belongs to $H^i(\calo_v,\calf)$ for all $v \in U- V$.
Then $\alpha$ is in the
image of the restriction map $H^i(U,\calf) \to H^i(V,\calf)$.
\end{lem}

\dem
The statement follows from the commutative diagram
$$
\begin{CD}
H^i(U,\calf ) @>>>  H^i(V,\calf) @>{\partial_v}>> \bigoplus_{v \in U-V} H^{i-1}(\kappa(v),F(-1)) \cr
@VVV @VVV @VV{=}V \cr
\bigoplus_{v \in U-V} H^i(\calo_v, \calf) @>>> \bigoplus_{v \in U-V} H^i(K_v,F) @>{\partial_v}>> \bigoplus_{v \in U-V} H^{i-1}(\kappa(v),F(-1))
\end{CD}
$$
whose exact rows come from localization sequences in \'etale cohomology.
\enddem

\noindent{\em Proof of Proposition \ref{partpoitou}.} Let $U_0$ and $\calf$ be as above. For
every nonempty Zariski open subset
$V \subset U_0$ there is an exact sequence of finite groups
$$
\cdots\to H^i_c(V,\calf)\to H^i(V, \calf)\to\bigoplus_{v\in X\setminus V} H^i(K_v, F)\to H^{i+1}_c(V, \calf)\to\cdots
$$
constructed as in (\cite{adt}, Lemma II.2.4), noting that we may replace henselisations at $v$ by completions in view of Greenberg's approximation theorem \cite{greenberg} (see e.g. the proof of \cite{dhsza1}, Lemma 2.7 for a detailed argument). Lemma~\ref{purlemma} then yields an exact sequence
\begin{equation}\label{varhensel}
H^i(U,\calf) \to \prod_{v \not \in U} H^i(K_v,F) \times
\prod_{v \in U-V} H^i(\calo_v,\calf) \to H^{i+1}_c(V,\calf)
\end{equation}
for each pair of nonempty open subsets $V\subset U$ contained in $U_0$. Moreover, the group $H^{i+1}_c(V,\calf)$ is dual to $H^{3-i}(V,\calf ')$ by Artin--Verdier duality for finite modules (see e.g. \cite{geisser}, \S 5.2). Therefore by taking the inverse
limit over $V$ we get an exact sequence
$$H^i(U,\calf) \to \prod_{v \not \in U} H^i(K_v,F) \times
\prod_{v \in U} H^i(\calo_v,\calf) \to H^{3-i}(K,F')^D.$$
(In fact, one needs to be a bit careful here since the index set of the inverse limit is uncountable. We have to check that if an element  $(f_v)$ of the middle term  has trivial image
in $H^{3-i}(K,F')^D$, then it comes from $H^i(U,\calf)$. For each nonempty open $V \subset U$
the assumption implies that the image of the truncated element
$(f_v)_{v \not \in V}$ in $H^{3-i}(V,\calf')^D$ is zero, so the inverse
image $E_V$ of $(f_v)_{v \not \in V}$ in $H^i(U,\calf)$ is
nonempty by exact sequence (\ref{varhensel}). On the other hand for $W \subset V$ we obviously have
$E_W \subset E_V$, so the intersection of all subsets $E_V$ in the finite set $H^i(U,\calf)$ is nonempty, as required.)

Finally, taking now direct limit over $U$ yields the exact sequence
$$H^i(K,F) \to \P^i(F) \to H^{3-i}(K,F')^D.$$
\enddem

The full statement of the Poitou--Tate sequence is:

\begin{theo}[Poitou--Tate sequence for finite modules] \label{poitfini}
Let
$F$ be a finite Galois module over $K$. There is a 12-term exact
sequence of topological groups:
$$
\begin{CD}
0 @>>> H^0(K,F) @>>> \prod_{v \in X^{(1)}} H^0(K_v,F) @>>> H^3(K,F')^D \cr
&&&&&& @VVV \cr
&& H^2(K,F') ^D @<<< \P^1(F) @<<< H^1(K,F) \cr
&& @VVV \cr
&& H^2(K,F) @>>> \P^2(F) @>>> H^1(K,F')^D \cr
&&&&&& @VVV \cr
0 @<<< H^0(K,F') ^D @<<< \bigoplus_{v \in X^{(1)}} H^3(K_v,F) @<<< H^3(K,F)
\end{CD}
$$
\end{theo}\smallskip

In this statement the groups $H^i(K,F)$ and $H^i(K,F')$
are equipped with the discrete topology (hence the dual of $H^i(K,F')$
is compact), and the topology of the restricted products has already been described.

For the proof we need the following general lemma about topological groups which is well known to the experts (whoever they may be) but for which we could not find a reference.

\begin{lem}\label{toplemma}
Let
$$0 \to A \stackrel{i}{\to} B \stackrel{p}{\to} C \to 0$$
be an exact sequence of  Hausdorff,
locally compact and totally disconnected topological abelian groups. Assume that
the morphisms $i$ and $p$ are strict. Then the dual sequence
$$0 \to C^D \to B^D \to A^D \to 0$$
is exact.
\end{lem}

Recall that a continuous homomorphism $f : A \to B$ of topological groups
is said to be {\it strict} if $f$ induces an isomorphism
between $A/\ker (f)$ (equipped with the quotient topology) and
$\im (f)$ (equipped with the subspace topology induced by $B$).
\smallskip

\begin{dem} The injectivity of the map $C^D \to B^D$ is obvious. Since $i$ is
strict, we can identify $A$ with its image in $B$, equipped
with the subspace topology.
A continuous homomorphism  $B \to \Q/\Z$ that is trivial on
 $A$ factors through the topological quotient $B/A$. As $p$ is strict, it thus  induces a continuous homomorphism $C \to \Q/\Z$. This shows the exactness of the dual sequence in the middle, so it remains to justify the surjectivity of
the map $B^D \to A^D$, which we now do following an argument kindly explained to us by
Yves de Cornulier.

Our task is to extend a continuous homomorphism $f : A \to \Q/\Z$ to a continuous homomorphism
$B \to \Q/\Z$. This is easy if $A$ is open in $B$: we may extend $f$ to a homomorphism $g: B \to \Q/\Z$ because $\Q/\Z$ is divisible,
and $g$ is automatically continuous because it is continuous in a neighbourhood of
$0$ by openness of $A$.

To reduce to the above special case, note first that by replacing $A$ and $B$ by their quotients modulo $\ker (f)$ we may assume
that $f$ is injective. This then implies that $A$ has the discrete topology because
$\ker (f)=\{0\}$ is then an open subgroup of $A$, the group $\Q/\Z$
being discrete. Now the assumption that $B$ is locally compact and
totally
disconnected implies by (\cite{hewitt}, Theorem 7.7)
that there is a basis of neighborhoods
of 0 in $B$ consisting of open subgroups. As $A$ is discrete,
we therefore find an open subgroup $U$ of $B$
such that $U \cap A=\{0\}$. Consider the union in $B$ of all open cosets $a+U$ for $a\in A$, and denote it by $A+U$. It is an open subgroup of $B$, and  there is a well-defined projection homomorphism $p_1:\, A+U\to A$ since $U\cap A=\{0\}$. The map $p_1$ is also continuous, the preimage of an open subset $V\subset A$ being the union of the open cosets $v+U$ for $v\in V$. Therefore the composite map $f\circ p_1:\, A+U\to \Q/\Z$ is a continuous homomorphism extending $f$. Replacing $A$ by $A+U$ we reduce to the open case treated above.
\end{dem}

\noindent{\em Proof of Theorem \ref{poitfini}.\/} Exactness of the two middle rows follows from Proposition \ref{partpoitou}. That  of the third is also a consequence, once we note the surjectivity of the map
$$\P^3(F)=\bigoplus_{v \in X^{(1)}} H^3(K_v,F) \to H^0(K,F')^D.$$
This follows from the injectivity of the dual maps
$H^0(K,F') \to  H^0(K_v,F')$ via the local duality
theorem for finite modules. The first row becomes the dual of the last one if we exchange $F$ and $F'$, so it is exact as well, the last row being an exact sequence of discrete torsion groups. The vertical maps
are defined using  Theorem~\ref{dualsha}.
It remains to show that the dual of the exact sequence
$$0 \to \Sha^i(F') \to H^i(K,F') \to \P^i(F')$$
is still exact for $i=1,2,3$. Since local duality (\ref{finitepairing}) induces a perfect duality of the locally compact groups $\P^i(F)$ and $\P^{3-i}(F')$ (this uses Proposition \ref{annihil} for $i=1,2$), we can then conclude the exactness of
$$\P^{3-i}(F) \to H^i(K,F')^D \to \Sha^{4-i}(F) \to 0.$$
To do so, we apply Lemma \ref{toplemma}, so we only have to prove that the morphism
$H^i(K,F') \to \P^i(F')$ is strict. It will suffice to show that its image is discrete.
For each open subset
$U \subset U_0$ (where $U_0$ is chosen such that
$F'$ extends to a finite and \'etale group scheme
$\calf'$ over $U_0$), the group $H^i(U,\calf')$ is finite. Since
each element of the image of $ H^i(K,F') $ in $\prod_{v \not \in U} H^i(K_v,F') \times
\prod_{v \in U} H^i(\calo_v, \calf')$ comes from $H^i(U,\calf')$ by
Lemma~\ref{purlemma}, we are done.\enddem

We now turn to tori. The following proposition is similar to
Proposition~\ref{partpoitou}.

\begin{prop}\label{before}
Let $T$ be a $K$-torus. There are exact sequences of topological groups
\begin{equation} \label{toreun}
H^1(K,T) \to \P^1(T) \to H^1(K,T')^D
\end{equation} and
\begin{equation} \label{toredeux}
H^2(K,T) \to \P^2(T)_{\tors} \to (H^0(K,T')_{\wedge})^D \to 0.
\end{equation}
\end{prop}\smallskip

Here $\P^1(T)$ is a topological  restricted product of finite discrete groups. The topology on $\P^2(T)_{\tors}$ is defined as follows. For each $n>0$ the group $_n \P^2(T)$ (restricted product of
the $_n H^2(K_v,T)$
with respect to the $_n H^2(\calo_v,\calt)$)
is equipped with the restricted product topology associated with the
discrete topology on each ${}_n H^2(K_v,T)$.
Their direct limit  $\P^2(T)_{\tors}$ is equipped with the
direct limit topology  and {\it not} by the subspace
topology induced by $\P^2(T)$.\medskip

Before starting the proof, we have to recall from \cite{pchp} some properties of the cohomology of tori over open subcurves of $X$.

\begin{lem} \label{firstprop} Let $U\subset X$ be a nonempty open subset such that $T$ extends to a torus $\T$ over $U$.

\noindent $(1)$ {\rm} For each $i \geq 1$ the group
$H^i(U,\T)_{\rm tors}$ is of cofinite type, and so is the compact support cohomology group
$H^i_c(U,\T)_{\rm tors}$ for $i \geq 2$.

\noindent$(2)$ For $U$ sufficiently small the group $H^2_c(U,\T)\{\ell\}$
is finite.

\noindent $(3)$ There is a canonical pairing
$$
H^i(U, \calt')\otimes H^{3-i}_c(U, \calt)\to\Q/\Z
$$
inducing a perfect pairing of finite groups
$$\overline{H^1(U,\T')\{\ell\}}  \times\overline {H^{2}_c(U,\T)}\{ \ell \}
 \to \Q_{\ell}/\Z_{\ell}$$
for each prime number $\ell$.

\end{lem}

\begin{dem}
The first statement is (\cite{pchp}, Proposition 3.4 (1)). To prove (2), note first that $H^2_c(U,\T)$ is a torsion group (\cite{pchp}, Corollary 3.3). Propositions 4.2 and 3.6 of \cite{pchp} imply that there is a well-defined map
$$
\bigoplus_{v \in X^{(1)}} H^1(K_v,T) \to H^2_c(U,\T)
$$
whose cokernel has finite $\ell$-primary part for $U$ sufficiently small. Since the groups $H^1(K_v,T)$ have bounded exponent by Hilbert's Theorem 90 and a restriction-corestriction argument, the finiteness of $H^2_c(U,\T)\{\ell\}$ follows from statement (1).
Finally, the last statement follows from the Artin--Verdier style result of (\cite{pchp}, Theorem 1.3), noting that for an $\ell$-primary torsion group $A$ of cofinite type the quotient $\overline A$ modulo the maximal divisible subgroup is isomorphic to the inverse limit $\varprojlim A/\ell^mA$.
\end{dem}

\noindent{\em Proof of Proposition \ref{before}.} Extend $T$, $T'$ to tori $\calt$, $\calt '$ over a nonempty
Zariski affine open subset $U_0 \subset X$.
We first prove the first statement along the lines of
Proposition~\ref{partpoitou}. The $K$-torus $T$ is split by an
extension of degree $m >0$.
Given nonempty open subsets $V \subset U \subset U_0$, there is an
exact sequence
$$H^1(U,\T) \to \prod_{v \not \in U} H^1(K_v,T) \times
\prod_{v \in U-V} H^1(\calo_v,\T) \to  H^2_c(V,\T)$$
obtained similarly as (\ref{varhensel}), except that we use Harder's lemma (\cite{harderinv}, Lemma 4.1.3) instead of Lemma \ref{purlemma}. The $K$-torus $T$ is split by some field
extension of degree $m >0$, so the middle term has exponent dividing $m$ by Hilbert's Theorem 90. We therefore obtain an exact sequence
$$H^1(U,\T)/m \to \prod_{v \not \in U} H^1(K_v,T) \times
\prod_{v \in U-V} H^1(\calo_v,\T) \to\,  _m H^2_c(V,\T)$$
whose terms are finite by Proposition ~\ref{firstprop} (1).
We
know by Lemma \ref{firstprop} (2)
that for $V$ sufficiently small the group $H^2_c(V,\T)\{\ell\}$
is finite, and therefore its maximal divisible subgroup is trivial. Lemma \ref{firstprop} (3) then implies
that the map $H^2_c(V,\T)\{\ell\} \to H^1(V,\T')\{\ell\} ^D$ is injective,
whence an exact sequence of finite groups
$$H^1(U,\T)/m \to \prod_{v \not \in U} H^1(K_v,T) \times
\prod_{v \in U-V} H^1(\calo_v,\T) \to \prod_{\ell \mid m}
H^1(V,\T')\{\ell\}^D.$$
Taking the inverse limit over $V \subset U$ as in the proof of Proposition \ref{partpoitou} and noting that $H^1(K,T')$ is also of exponent dividing $m$, we get an exact
sequence
$$H^1(U,\T)/m \to \prod_{v \not \in U} H^1(K_v,T) \times
\prod_{v \in U} H^1(\calo_v,\T) \to H^1(K,T')^D.$$
Replacing $H^1(U,\T)/m$ by $H^1(U,\T)$ and taking the direct
limit along the subsets $U$ we then obtain exact sequence  (\ref{toreun}).

To prove the second statement, we fix an integer $n >0$ and start with the exact sequence
$$H^2(K,\,{}_n T) \to \P^2(_n T) \to H^1(K,\,{}_n T')^D \to H^3(K,\,{}_n T)$$
given by Theorem~\ref{poitfini}. Consider the commutative diagram with exact first row
$$
\begin{CD}
\varinjlim H^2(K,_n T) @>>> \varinjlim \P^2(_n T) @>>> \varinjlim
(H^1(K,_n T')^D)  \cr
@VVV @VVV @VVV \cr
H^2(K,T) @>>> \P^2(T)_{\tors} @>>> (H^0(K,T')_{\wedge})^D
\end{CD}
$$
where the dual of the third vertical map comes from the exact sequence
$$0 \to H^0(K,T')/n \to H^1(K,_n T') \to {}_n H^1(K,T') \to 0.$$
The map in question is in fact an isomorphism because the Tate module of $H^1(K,T')$ is zero,
the group $H^1(K,T')$ being of finite exponent. Similarly, the first vertical map is an isomorphism
because $H^{1}(K,T)$ is torsion and hence $H^{1}(K,T) \otimes \Q/\Z=0$. Since the maps $H^2(K_v,_n T) \to {}_nH^2(K_v,T)$ and
$H^2(\calo_v, _n \calt) \to {}_nH^2(\calo_v,\calt)$ are onto
(the latter for $v \in U_0$), the middle vertical map is also onto. A diagram chase then yields exact sequence (\ref{toredeux}) except for the zero on the right. By the diagram and by the exact sequence of Theorem \ref{poitfini} this will follow if we prove that the direct limit over $n$ of the groups $H^3(K,\, {}_nT)$ is trivial. Note first that this direct limit is $H^3(K,T)$ because $H^{2}(K,T) \otimes \Q/\Z=0$. But $H^3(K,T)$ is the direct limit of the groups $H^3(V,\calt)$ for $V\subset U_0$ and the latter groups are trivial by
(\cite{schvh}, Corollary~4.10).\enddem

We shall also need the following consequence of local duality.

\begin{lem} \label{restrtore}
The local duality theorem induces a perfect duality between
the locally compact groups $\P^1(T)$ and $\P^1(T')$. It also
induces an isomorphism $$(\P^2(T)_{\tors}) ^D \simeq \P^0(T')_{\wedge}. $$
\end{lem}

\dem The first assertion is an immediate consequence of
Prop.~\ref{annihiltori} $a)$. The second one follows from part $b)$ of the proposition, the definition of the topology on $\P^2(T)_{\tors}$ as well as
the canonical isomorphism between $\P^0(T)/n$
and the restricted product of the $H^0(K_v,T)/n$ for each $n>0$. To see that the natural surjections $\P^0(T)/n\to H^0(K_v,T)/n\/$ indeed induce an isomorphism with the restricted product it suffices to show the injectivity of the maps $H^0(\calo_v,\calt)/n\to H^0(K_v,T)/n$. This in turn follows from the injectivity of the maps $H^1(\calo_v, {}_n\calt)\to H^1(K_v, {}_nT)$ in view of the Kummer sequence.
\enddem

Finally, we shall use an analogue of Lemma \ref{purlemma} (again valid over an arbitrary base field).

\begin{lem} \label{harder2}
Assume $T$ extends to a torus $\calt$
over a nonempty
Zariski open subset $U\subset X$, let $V\subset U$ be a nonempty open subset, and $i\geq 2$ an integer. If $\alpha \in H^i(V,T)$ is such that its restriction
${\alpha_v\in H^i(K_v,T)}$ lies in $H^i(\calo_v,\calt)$ for each $v \in U-V$, then $\alpha$ is
in the image of $H^i(U,\calt)$.
\end{lem}

\begin{dem} Lift $\alpha$ to some $\alpha_n \in H^i(V,\, _n T)$ for $n >0$.
For $v \in U-V$ the image of the restriction $\alpha_{n,v} \in H^i(K_v,\, _n T)$
in $H^i(K_v,T)$ comes from an element in $ H^i(\calo_v,\calt)$.
We have a commutative diagram with injective vertical maps:
$$
\begin{CD}
\varinjlim_n \, H^i(\calo_v,\, _n \calt) @>{\simeq}>> H^i(\calo_v, \calt) \cr
@VVV @VVV \cr
\varinjlim_n \, H^i(K_v,\, _n T) @>{\simeq}>> H^i(K_v, T).
\end{CD}
$$
The horizontal maps are isomorphisms because $H^{i-1}(\calo_v, \calt)\otimes\Q/\Z=0$ for $i\geq 2$.
Since there are only finitely many places $v$ in $U-V$,
there exists a multiple $m$ of $n$ such that
the image $\alpha_m$ of $\alpha_n$ in $H^i(V,_m \calt)$
has the following property: for each $v \in U-V$, the restriction
$\alpha_{m,v}$ belongs to $H^i(\calo_v,_m \calt)$. By Lemma~\ref{purlemma}
this implies
that $\alpha_m$ is in the image of $H^i(U,_m \calt)$, hence $\alpha$ is
in the image of $H^i(U,\calt)$.
\end{dem}

We now arrive at

\begin{theo}[Poitou--Tate sequence for tori] \label{poitores}
Let
$T$ be a $K$-torus. There is an exact
sequence of topological groups:
$$
\begin{CD}
0 @>>> H^0(K,T)_{\wedge}  @>>> \P^0(T)_{\wedge}  @>>> H^2(K,T')^D \cr
&&&&&& @VVV \cr
&& H^1(K,T') ^D @<<< \P^1(T) @<<< H^1(K,T) \cr
&& @VVV \cr
&& H^2(K,T) @>>> \P^2(T)_{\tors} @>>> (H^0(K,T')_{\wedge})^D @>>> 0
\end{CD}
$$
\end{theo}

\smallskip

\begin{dem} Exactness of the second and third rows follows
from Proposition \ref{before}. The vertical maps are defined using
Theorem~\ref{dualsha}. The exactness of
$$\P^1(T) \to H^1(K,T')^D \to \Sha^2(T) \to 0 $$
is proven as in
Theorem~\ref{poitfini}, replacing the use of Proposition \ref{annihil} by that of Proposition \ref{annihiltori} and
the finiteness of $H^1(U,\calf')$ by the finiteness of $H^1(U,\calt')/m$,
where $T'$ is split by an extension of degree $m$.

To get the exactness of
$$\P^0(T)_{\wedge} \to H^2(K,T')^D \to \Sha^1(T) \to 0 ,$$
we observe that by Lemma~\ref{restrtore}
it is the dual of the exact sequence
$$0 \to \Sha^2(T') \to H^2(K,T') \to \P^2(T')_{\tors}.$$
To see that the dual remains exact, we need to check, again using Lemma \ref{toplemma}, that $H^2(K,T')\to\P^2(T')_{\tors}$ is a strict morphism. Again we show that in fact the image $I$ of $H^2(K,T')$  is discrete in $\P^2(T')_{\tors}$. Fix $n >0$, and
let $A_n$ be the inverse image in $H^2(K,T')$
of the subgroup $B_n$ of $\P^2(T')_{\tors}$ defined by
$$B_n:= \prod_{v \not \in U} \, _n H^2(K_v,T') \times \prod_{v \in U} \,
_n H^2(\calo_v, \calt'). $$
By Lemma~\ref{harder2}, the group $A_n$ lies in the image of the map $H^2(U,\calt')
\to H^2(K,T')$ (lift first elements of $A_n$ to $H^2(V,\calt')$ for $V$ sufficiently small). Therefore $A_n$ is a torsion group of cofinite type by Proposition \ref{firstprop} (1),
which implies that $A_n/n$ is finite. In particular the image of $A_n$
in the $n$-torsion group $B_n$ is finite. Since this image is  $I \cap B_n$, we obtain that $I$ is discrete in $\P^2(T')_{\tors}$ by definition of the topology
on $\P^2(T')_{\tors}$.

Finally, the
exactness of the first row is obtained by dualizing the last one
(after exchanging $T$ and $T'$), once we have observed that
for each $n>0$ the bidual of the discrete torsion group
$H^0(K,T')/n$ is itself. The dual of the last row remains exact
because by definition of
the topology on $\P^2(T)_{\tors}$ and the equality  $(H^0(K,T')_{\wedge}) ^D=\varinjlim_n \,
(H^0(K,T')/n)^D$ the surjection $\P^2(T)_{\tors} \to (H^0(K,T')_{\wedge}) ^D$
is strict.
\end{dem}

\begin{rema}
{\rm The above sequence looks like the Poitou-Tate exact sequence for a finite
module (or a torus)
over a number field. But there is an important difference: the
group $\P^2(T)_{\tors}$ is not discrete (neither is its quotient
by the image of $H^2(K,T)$), and  its dual is not compact.
}
\end{rema}

\section{Obstruction to weak approximation for tori}\label{wa}

In this section we address weak approximation questions for tori.
First some notation: for a finite set of places $S \subset X^{(1)}$ and a positive
integer $i$ we define
$$\Sha^i_S(F):=\ker [H^i(K,F) \to \prod_{v \not \in S} H^i(K_v,F)]$$
for every finite group scheme $F$ over $K$.
We also define $\Sha^i_{\omega}(F)\subset H^i(K,F)$ as the subgroup of
 elements whose restriction to
$H^i(K_v,F)$ is zero for all but finitely many $v \in X^{(1)}$. Thus $\Sha^i_{\omega}(F)$
is the direct limit of the $\Sha^i_S(F)$ over all finite $S \subset X^{(1)}$. Similar notation
will be used for tori.

We start with the following consequence of Proposition~\ref{partpoitou}.

\begin{lem} \label{poitouS}

Let $S \subset X^{(1)}$ be a finite set of places, and let $F$ be a finite group scheme over $K$.
For $i=1,2$ there are exact sequences
$$H^i(K,F) \to \prod_{v \in S} H^i(K_v,F) \to \Sha^{3-i}_S(F')^D$$ where the second map is
given by
$$(f_v) \mapsto (g \mapsto \sum_{v \in S} \langle f_v, g_v \rangle_v ) \quad\text{for}\quad
(f_v) \in \prod_{v \in S} H^i(K_v,F),\, g \in \Sha^{3-i}_S(F'). $$

\end{lem}

\begin{dem}  Proposition~\ref{partpoitou} applied to $F'$ implies an exact sequence
$$\Sha^i_S(F') \to \prod_{v \in S} H^i(K_v,F') \to H^{3-i}(K,F)^D$$
because an element in the middle group can be completed by 0's to yield an element of ${\bf P}^i(F')$.
 Now the result is obtained by dualizing, replacing $i$ by $3-i$ and applying local duality  (\ref{finitepairing}).
\end{dem}

Next comes the key finiteness lemma. Its first part generalizes Remark 3.5 of \cite{pchp}.

\begin{lem} \label{finishaS}
Let $T$ be a $K$-torus, and $S \subset X^{(1)}$ a finite set of places.

a) If $T$ is quasi-trivial, then $\Sha^2_S(T)=\Sha^2_{\omega}(T)=0$.

\smallskip

b) In general the group  $\Sha^2_{\omega}(T)$ is of finite exponent and
$\Sha^2_S(T)$ is finite.
\end{lem}

\dem $a)$ (after J-L. Colliot-Th\'el\`ene).
Since $\Sha^2_{\omega}(T)$ is the union of all $\Sha^2_S(T)$ for $S$
finite, we only have to prove the statement about $\Sha^2_S(T)$.
By Shapiro's lemma it
is sufficient to deal with the case $T=\G$.
Each element $\alpha\in\Sha^2_S(\G)$
comes from an element
$\alpha_U\in\br U$ for some open $U\subset X$ whose restriction to $\br K_v$ is 0 for all $v \in U$. It is
sufficient to prove that $\alpha_U \in \br X$. Indeed, this will then imply $\alpha_U\in\br\calo_{X,v}$ for all $v\in X^{(1)}$, where $\calo_{X,v}$ is the local ring of $X$ at the closed point $v$ whose completion is the ring of integers $\calo_v$ of $K_v$. Since the  restriction map $\br\calo_v\to\br K_v$ is injective, we conclude that $\alpha_U$ maps to 0 already in $\br\calo_v$. Composing with the isomorphism $\br\calo_v\cong\br\kappa(v)$ we conclude that $\alpha_U$ maps to 0 by the `evaluation at $v$' map $\br X\to\br\kappa(v)$ induced by $v$.
Therefore $\alpha_U$ is orthogonal to $\pic X$ for Lichtenbaum's duality
pairing $\br X\times \pic X\to\Q/\Z$ which is defined in \cite{licht} by a sum of evaluation maps, from which we conclude $\alpha_U=0$.

Pick $v \in X-U$ and set $k':=\kappa(v)$.
We have to show
that the residue $\chi\in H^1(k',\Q/\Z)$
of $\alpha$ at $v$ is zero. Fix a positive integer $m$ such that
$\chi \in H^1(k',\Z/m)$.
Since $H^1(k',\Z/m)\cong H^1(\calo_{X,v}^h,\Z/m)$, we may lift $\chi \in H^1(k',\Z/m)$  to
an element $\tilde \chi \in H^1(\widetilde V,\Z/m)$
for an \'etale neighbourhood
$f : \widetilde V \to X$ of $v$
containing a point $w$ with
$f(w)=v$ and $\kappa(w)=k'$.
We may further assume $\widetilde V$ to be integral and write $\widetilde K$ for its function
field. Let $\pi$ be a uniformizing parameter of
the local ring $\calo_{\widetilde V,w}$.
We define an element $\beta$ of $\br \widetilde K$ by
$$\beta=f^*(\alpha)-(\pi \cup \tilde \chi) , $$
where $\pi$ is viewed in ${\widetilde K}^\times/{\widetilde K}^{\times m}=
H^1(\widetilde K,\mu_m)$.
By a standard residue computation, the residue of $\beta$ at
$w$ is zero, and  after shrinking $\widetilde V$ if necessary we
can assume that $\beta \in \br \widetilde V$.
Since $f$ is \'etale, by the implicit function theorem it induces
a homeomorphism
for the $p$-adic topology from an open neighborhood
$\widetilde \Omega$
of $w$ onto an open neighborhood $\Omega$ of $v$. We may further assume
$\Omega -\{v\}\subset U(k')$. Since
$\beta \in \br \widetilde V$,  we can assume
 (shrinking $
\widetilde\Omega$ if necessary) that
the evaluation map
$$\widetilde \Omega \to \br k' , \quad \widetilde P \mapsto
\beta(\widetilde P)$$
is constant (see for example \cite{bordem}, Lemma~6.2).
Then, using the assumption $\alpha_{U,u}=0$ for $u \in U$,
we see that
the evaluation map
$$\widetilde \Omega-\{ w \} \to \br k' , \quad \widetilde P
\mapsto f^*(\alpha)(\widetilde P)$$
is identically 0, so the map
$$\widetilde \Omega-\{ w \} \to \br k' , \quad \widetilde P \mapsto
(\pi \cup \tilde \chi)(\widetilde P)$$
is constant.
\smallskip

Shrinking again  $\widetilde V$ if necessary, we can find a Zariski open
subset $W\subset \A^1_{k'}=\spec (k'[t])$ containing $0$  and a smooth
$k'$-morphism $\varphi : \widetilde V \to W$ sending $w$ to 0 such that the pullback of $t\in \calo_{W,0}$ to $\calo_{\widetilde V,w}$ is $\pi$.
Using again the implicit function theorem, we find a $p$-adic neighborhood
$\Omega_1$ of $0 \in W(k')$ such that the evaluation map
$$\Omega_1 - \{ 0\} \to \br k' , \quad P \mapsto (t \cup \chi)(P)=P \cup \chi$$
is constant. In fact, it is identically 0 because if we take a uniformizing parameter $\pi_{k'}\in\calo_{k'}$ and consider $P=\pi_{k'}^{rm}$ for $r$ large enough, then $P\in\Omega_1$ and also $P\cup\chi=0$ as $m\chi=0$.
Now considering elements of the form $u.\pi_{k'}^n$ for large $n$ with some unit $u$, we see that we may represent classes in $H^1(k',\mu_m)={k'}^\times/
{k'}^{\times m}$ by points in $\Omega_1$, and hence we conclude that  $\chi$ is orthogonal to $H^1(k',\mu_m)$ for Tate's local duality pairing. This implies $\chi=0$, as desired.

\smallskip

$b)$ We deduce via a restriction-corestriction argument from the case
$T=\G$ that the groups $\Sha^2_S(T)$ and
$\Sha^2_{\omega}(T)$ are of finite exponent. On the other hand,
using the exact sequence
$$0 \to \Sha^2(T) \to \Sha^2_S(T) \to \bigoplus_{v \in S}
H^2(K_v,T) ,$$
we obtain that $\Sha^2_S(T)$ is of cofinite type because
$\Sha^2(T)$ is finite by Proposition~\ref{firstprop} (2)
and for each $n >0$ the Kummer sequence induces a surjection of the finite group  $H^1(K_v,_n T)$ onto ${}_n H^2(K_v,T)$.
\enddem

 We now turn to statements concerning weak approximation. Given a variety $Y$ over $K$, we equip the sets $Y(K_v)$ with the $v$-adic topology defined by the discrete valuation coming from the closed point $v\in X$. By definition, weak approximation holds for $Y$ if the diagonal image of $Y(K)$ is dense in the topological direct product of the $Y(K_v)$. Note that this is the case when $Y$ is smooth and $K$-rational. Indeed, for the affine space weak approximation is just the classical Artin--Whaples approximation lemma (\cite{lang}, Theorem XII.1.2) for discrete
valuations, and the general case follows from the implicit function theorem for ultrametric valuations
(\cite{serrelie}, Part II, \S III.10) which implies
that for all $v$ and all nonempty Zariski open subsets $U\subset X$ the subset $U(K_v)$
is dense in $X(K_v)$ for the $v$-adic topology.

\begin{theo} \label{mainwa}
Let $T$ be a torus over $K$.

\smallskip

a) For a finite set $S \subset X^{(1)}$  of places
denote by $\overline{T(K)}_S$ be the closure of $T(K)$ in
$\prod_{v \in S} T(K_v)$.
There is an exact sequence of finite groups
\begin{equation} \label{sansvoskS}
0 \to \overline{T(K)}_S \to \prod_{v \in S} T(K_v) \stackrel{\lambda_S}\longrightarrow \Sha^2_S(T')^D
\to \Sha^1(T) \to 0 .
\end{equation}
where the map $\lambda_S$ is
defined using the local duality pairings $\langle \,\,, \,\, \rangle _v$ by
$$(t_v) \mapsto (t' \mapsto \sum_{v \in S} \langle t_v, t'_v \rangle _v)$$
for $t_v \in T(K_v)$ and $t' \in \Sha^2_S(T').$\smallskip

b) Similarly, denote by $\overline{T(K)}$ the closure of
$T(K)$ in the (topological) product of the $T(K_v)$ for all $v$. There is an exact sequence

\begin{equation} \label{sansvosk}
0 \to \overline{T(K)} \to \prod_{v \in X^{(1)}} T(K_v)\stackrel{\lambda_\omega} \longrightarrow
\Sha^2_{\omega}(T')^D \to \Sha^1(T) \to 0
\end{equation}
with an analogously defined map $\lambda_\omega$.
\end{theo}

Note that the maps $\lambda_S$ and $\lambda_\omega$ are continuous because the local duality pairings $H^0(K_v,T)\to H^2(K_v, T')^D$ factor through the profinite completions $H^0(K_v,T)^\wedge$. \medskip

\dem a) We follow the method of Sansuc \cite{sansuc}. First, note that
a quasi-trivial $K$-torus satisfies weak approximation as it is $K$-rational, and therefore the cokernel of the first map in (\ref{sansvoskS}) does not change if we multiply $T$ by a quasi-trivial torus. On the other hand, Lemma \ref{finishaS} $a)$ and Hilbert's Theorem 90 imply that the last two terms in the sequence do not change either when  $T$ is multiplied by a quasi-trivial torus. So we are free to perform such an operation during the proof. Also, if the statement holds for some finite direct power $T^m$ of $T$, it holds for $T$ as well. So up to replacing $T$ by some $T^m\times Q$ with $Q$ quasi-trivial, we may apply a lemma of Ono  (see e.g.\cite{sansuc}, Lemma~1.7) to find
an exact sequence of commutative $K$-group schemes
$$0 \to F \to R \to T \to 0$$
where $R$ is a quasi-trivial torus and $F$ is finite. This induces
a dual exact sequence of Galois modules
$$0 \to F' \to T' \to R' \to 0 $$
where $F'=\Hom(F, \Q/\Z(2))$ as usual. (Indeed, the map $T'\to R'$ is surjective since the map of cocharacters $\widecheck R\to\widecheck T$ is injective by finiteness of $F$. On the other hand, tensoring the exact sequence of character groups
$$0  \to\widehat T \to\widehat R \to \Hom(F,\Q/\Z(1))\to 0 $$
by $\Z(1)[1]$ in the derived sense we obtain an exact triangle
$$
\Hom(F,\Q/\Z(2))\to \widehat T(1)[1]\to \widehat R(1)[1]\to \Hom(F,\Q/\Z(2))[1]
$$
in which using formula (\ref{dualtordef}) we may identify the middle map with the surjection $T'\to R'$.)
There is a commutative diagram with exact rows~:
\begin{equation} \label{diagweak}
\begin{CD}
R(K) @>>> T(K) @>>> H^1(K,F) @>>> 0 \cr
@V{\rho_R}VV @VVV @VVV \cr
\prod_{v \in S} R(K_v) @>>> \prod_{v \in S} T(K_v) @>>> \prod_{v \in S}
H^1(K_v,F) @>>> 0 \cr
&& @VVV @VVV \cr
&& \Sha^2_S(T')^D @>\simeq>> \Sha^2_S(F')^D
\end{CD}
\end{equation}
where the isomorphism at the bottom comes from the vanishing of the groups
$H^1(K,R')$ and $ \Sha^2_S(R')$ for the quasi-trivial torus $R'$ (Lemma \ref{finishaS} $a)$).

\smallskip

The last column is exact by Lemma~\ref{poitouS}. Moreover, the map $\rho_R$ has dense image, again
 because $R$ is a quasi-trivial
torus.
A diagram chasing therefore yields an exact sequence
$$0 \to \overline{T(K)}_S \to \prod_{v \in S} T(K_v) \to \Sha^2_S(T')^D .$$

\smallskip

To get the exactness of the last three terms in (\ref{sansvoskS}), we start
from the exact sequence
$$0 \to \Sha^2(T') \to \Sha^2_S(T') \to \bigoplus_{v \in S}
H^2(K_v,T') .$$
Dualizing this sequence, we get using Proposition~\ref{localdual}
and Theorem~\ref{dualsha} an exact sequence
$$\prod_{v \in S} T(K_v)^{\wedge} \to \Sha^2_S(T')^D \to
\Sha^1(T) \to 0.$$
But $\prod_{v \in S} T(K_v)$ and its profinite completion
$\prod_{v \in S} T(K_v)^{\wedge}$
have the same image in
$\Sha^2_S(T')^D$, because the latter is a finite group by
Lemma~\ref{finishaS} $b)$.

\smallskip

b) Taking the projective limit over $S$ in (\ref{sansvoskS}),
we obtain the exactness of
$$0 \to \overline{T(K)} \to \prod_{v \in X^{(1)}} T(K_v)
\to \Sha^2_{\omega}(T')^D .$$
To get the exactness of the last three terms, we proceed as in $a)$.
After dualizing the exact sequence
$$0 \to \Sha^2(T') \to \Sha^2_{\omega}(T') \to \bigoplus_{v \in X^{(1)}}
H^2(K_v,T') ,$$
we see that it is sufficient to prove that the image of
$\prod_{v \in X^{(1)}} T(K_v)$ in $\Sha^2_{\omega}(T')^D$ is closed.
It is therefore sufficient to show that the quotient of
$\prod_{v \in X^{(1)}} T(K_v)$ by $\overline{T(K)}$ is compact.
Using diagram (\ref{diagweak}), we obtain that this quotient is the same as
the quotient of the profinite group $\prod_{v \in X^{(1)}} H^1(K_v,F)$
by the closure of the image of $H^1(K,F)$, hence it is compact as well.
This concludes the proof.
\enddem

Recall that $T$ is said to satisfy the {\em weak weak approximation
property} if there exists a finite set $S_0 \subset X^{(1)}$ such
that for every finite set $S \subset X^{(1)}$ with $S \cap S_0=\emptyset$
the set $T(K)$ is dense in $\prod_{v \in S} T(K_v)$.

\begin{cor}
A $K$-torus $T$ satisfies the weak weak approximation
property  if and only if
$\Sha^2_{\omega}(T')$ is finite.
\end{cor}

\dem Assume first $\Sha^2_{\omega}(T')$ finite, and let $S_0$ be the union of all
places $v$ such that there exists $\alpha \in \Sha^2_{\omega}(T')$ with
$\alpha_v \neq 0$. Then $S_0$ is finite. If $S \subset X^{(1)}$
is a finite set with $S \cap S_0=\emptyset$, then
by construction we have
$\Sha^2_{S}(T')=\Sha^2(T')$, which implies that $T(K)$ is dense in
$\prod_{v \in S} T(K_v)$ by Theorem \ref{mainwa} $a)$ and Theorem \ref{dualsha}.

Conversely, if there exists a finite set $S_0 \subset X^{(1)}$ such
that $T$ satisfies weak approximation outside $S_0$, then for
every finite set of places $S$ disjoint from $S_0$ we have $\Sha^2_{S}(T')=\Sha^2(T')$ by Theorem \ref{mainwa} $a)$ and Theorem \ref{dualsha}. Hence there is an exact sequence
$$0 \to \Sha^2(T') \to \Sha^2_{\omega}(T') \to \bigoplus_{v \in S_0}
H^2(K_v,T') .$$
Now the finiteness of
$\Sha^2_{\omega}(T')$ follows in the same way as in the proof of Lemma~\ref{finishaS} $b)$.
\enddem

We close this section by observing that, in contrast to the number field case, it may happen that for a $K$-torus $T$
the group $\Sha^2_{\omega}(T)$ is infinite, and therefore the obstruction in Theorem \ref{mainwa} is not finite.

\begin{prop}\label{exsha}
The exists a torus $T$ over $K=\Q_p(t)$ with $\Sha^2_{\omega}(T)$ infinite.
\end{prop}

\begin{dem}
It suffices to find a torus $Q$ over $K=\Q_p(t)$ with infinite $\Sha^1_{\omega}(Q)$. Indeed, since for a quasi-trivial torus $P$ we have $H^1(K,P)=
\Sha^2_{\omega}(P)=0$, choosing a flasque resolution $1\to T\to P\to Q\to 1$ yields $\Sha^1_{\omega}(Q)\cong \Sha^2_{\omega}(T)$.

One can then choose $Q$ to be the base change to $K$ of a coflasque $\Q_p$-torus $Q_0$. Since $Q_0$ is coflasque, we
have  by definition $H^1(k',\widehat Q_0)=0$ for every finite extension $k'|\Q_p$, whence $H^1(k',Q_0)=0$ by Tate local duality for tori. This implies
$H^1(K,Q)=\Sha^1_{\omega}(Q)$. Indeed, each class in $\alpha\in H^1(K, Q)$ comes from a class in $H^1(U, {\mathcal Q})$ for a suitable open subset $U\subset {\bf P}^1_{\Q_p}$ and model $\mathcal Q$ over $U$, and thus for $v\in U$ the image of $\alpha$ in $H^1(K_v, Q)$ comes from $H^1(\calo_v, {\mathcal Q})\cong H^1(k(v), Q_0)$ which is 0 by the above.

Choose moreover $Q_0$ in such a way that
$H^1(\Q_p,\widehat Q_0') \neq 0$ for the dual flasque torus $Q_0'$.
For an explicit example, one may take $Q_0'$ to be the flasque torus in a
flasque resolution of
$\R^1 _{k'|\Q_p} \G$, with $G:=\gal(k'|\Q_p)$ isomorphic to $(\Z/2\Z)^2$.
Indeed, we have  $H^1(\Q_p,\widehat Q_0')\cong H^3(G,\Z)$ by (\cite{requ}, Proposition 7, page 191; see also Lemma 5 and Proposition 6 of the same paper), but
$H^3(G,\Z)\cong\Z/2\Z$  for our $G$ by (\cite{cogal}, I.4.4, Proposition 28).

Now fix a nonzero element $a \in H^1(\Q_p,\widehat Q_0') \neq 0$. For each
$b \in \Q_p$ consider the cup-product $A_b:={i(a)\cup (t-b)}$ in $H^1(K,Q)$, where
$(t-b)$ is viewed as an element of $H^0(K,\G)$ and $i$ is the composite map $H^1(k(b), \widehat Q')\stackrel{\rm inf}\to H^1(K_b, \widehat Q')\to  H^1(K, \widehat Q')$ (observe that $\widehat Q' \otimes
\G=Q$). The image of $A_b$  by the residue map $H^1(K,Q)\to H^1(k(b),\widehat Q')$ is $a$, so $A_b$ does not
come from $H^1(U_b,{\mathcal Q})$ for any open neighbourhood $U_b$ of $b$ (see the remark below for the residue argument used here). Since $\Q_p$-rational points are Zariski dense in $\P^1_{\Q_p}$, there is
no Zariski open subset $U\subset{\bf P}^1_{\Q_p}$ such that every element
of $H^1(K,Q)$ comes from $H^1(U,{\mathcal Q})$. Therefore $H^1(K,Q)$ must be infinite.
\end{dem}

\begin{rema}\rm Some explanation is in order concerning the residue maps used in the above proof. They appear in various guises in the literature, but the following simple direct approach which is well adapted for our purposes is perhaps worth noting. Given a complete discrete valuation ring $\calo_v$ with closed point $v$, fraction field $K_v$ and perfect residue field $k(v)$, it is shown in \cite{adt}, II.1.7(b) that there are isomorphisms $H^{i+1}_v(\calo_v,\G)\stackrel\sim\to H^{i}(k(v),\Z)$ for $i\geq 1$, basically induced by the valuation map. Given an $\calo_v$-torus $\T$, tensoring $\G$ with $\widecheck \T=\widehat \T'$ we may repeat the construction in the reference to obtain an isomorphism $H^{i+1}_v(\calo_v,T)\stackrel\sim\to H^{i}(k(v),\widecheck T)$. On the other hand, localization on $\spec\calo_v$ yields an exact sequence
$$
H^i(\calo_v,\T)\to H^i(K_v,T)\to H^{i+1}_v(\calo_v,\T)
$$
whence the residue map $H^i(K_v,T)\to H^{i}(k(v),\widecheck T)$ used above for $i=1$. By construction, for $i=1$ its value is $a$ on elements of the form ${\rm inf}(a)\cup \pi$ for $a\in H^1(k(v),\widecheck T)$ and $\pi$ a uniformizer, and it vanishes on classes coming from $H^1(\calo_v,\T)$.

In the above proof this construction was applied to completions of local rings of smooth $k$-curves. In fact, it can be shown using the global localization sequence in \'etale cohomology that on a smooth $k$-curve $U$ with fraction field $K$ the classes in $H^i(K,T)$ that come from some class in $H^i(U,\mathcal T)$ are exactly those whose residues are trivial at all closed points of $U$, but we did not need this fact.
\end{rema}

\section{Reciprocity obstructions to weak approximation}\label{recwa}

For varieties over number fields, unramified elements in the Brauer group have been used for a long time to construct obstructions to weak approximation. As was observed by Colliot-Th\'el\`ene in the 1990's, such obstructions can be defined using higher unramified cohomology groups as well. A detailed exposition is given in \S 2 of the recent preprint \cite{ctps2}; we give here a reasonably self-contained exposition of the special case we need.

Recall that given a field $F$ and a smooth $F$-variety $Y$ with function field $F(Y)$, the unramified part $H^i_{\rm nr}(F(Y), \mu_n^{\otimes j})$ of the cohomology group $H^i(F(Y), \mu_n^{\otimes j})$ is defined as the group of cohomology classes that come from $H^i(A, \mu_n^{\otimes j})$ for every discrete valuation ring $A\supset F$ with fraction field $F(Y)$. If moreover $Y$ is proper over $F$, the Bloch--Ogus theorem for the local rings of $Y$ implies that $H^i_{\rm nr}(F(Y), \mu_n^{\otimes j})$ consists of the classes that come from $H^i(\calo_{Y,P}, \mu_n^{\otimes j})$ for every point $P\in Y$ (see \cite{ctsb}, Theorem 4.1.1). As pointed out on p. 156 of \cite{ctcrelle}, for an overfield $F'\supset F$ this enables one to define an evaluation pairing
$$
Y(F')\times H^i_{\rm nr}(F(Y), \mu_n^{\otimes j})\to H^i(F', \mu_n^{\otimes j})
$$
as follows: given a point $M:\, \spec F'\to Y$ whose image is the point $P\in Y$, one evaluates $\alpha\in H^i_{\rm nr}(F(Y), \mu_n^{\otimes j})$ at $M$ by lifting it (uniquely) to $H^i(\calo_{X,P}, \mu_n^{\otimes j})$ and then taking its image by the map $H^i(\calo_{X,P}, \mu_n^{\otimes j})\to H^i(F', \mu_n^{\otimes j})$ induced by $M$. Note that this definition also works for $Y$ not necessarily proper but having a smooth compactification over $F$.

If now $F=K$ is the function field of a curve $X$ over a finite extension $k|\Q_p$ and $F'=K_v$ is the completion of $K$ at a closed point $v\in X$, for $i=3$ and $j=2$ the above construction yields a pairing
\begin{equation}\label{evalpairing}
Y(K_v)\times H^3_{\rm nr}(K(Y), \mu_n^{\otimes 2})\to H^3(K_v, \mu_n^{\otimes 2}).
\end{equation}
>From this we can construct a pairing
\begin{equation}\label{nrpairing}
\prod_{v \in X^{(1)}} Y(K_v)\times H^3_{\rm nr}(K(Y), \Q/\Z(2))\to \Q/\Z
\end{equation}
using that $ H^3(K_v, \Q/\Z(2))\cong \Q/\Z$ for all $v$, as recalled in Section \ref{one}.
To see that it is well defined, we show:

\begin{lem}\label{lemeval} For fixed $\alpha\in H^3_{\rm nr}(K(Y), \mu_n^{\otimes 2})$ the map $Y(K_v)\to\Q/\Z$ induced by evaluating $\alpha$ via $(\ref{evalpairing})$ is 0 for all but finitely many $v$.
\end{lem}

\begin{dem} This is part of (\cite{ctps2}, Proposition 2.5). The argument is the following. Choosing a model ${\mathcal Y}^c\to U$ of a smooth compactification $Y^c$ over a sufficiently small open $U\subset X$, the element $\alpha$ lies in $H^3(\calo_{{\mathcal Y}^c,Q}, \Q/\Z(2))$ for all but finitely many codimension 1 points $Q$ of the $k$-variety ${\mathcal Y}^c$. By the assumption on $\alpha$, the exceptional $Q$ lie in finitely many closed fibres ${\mathcal Y}^c_{v_1},\dots,{\mathcal Y}^c_{v_r}$ of ${\mathcal Y}^c\to U$. Therefore for $v\neq v_i$ we have $\alpha\in H^3(\calo_{{\mathcal Y}^c,P}, \Q/\Z(2))$ for all  points $P\in {\mathcal Y}^c_{v}$ and hence the map $Y^c(K_v)={\mathcal Y}^c(\calo_v)\to \Z/n\Z$ induced by $\alpha$ factors through $H^3(\calo_v, \Q/\Z(2))=0$.\end{dem}

By Saito (\cite{saito}, Proof of Proposition II.1.2), the sequence
$$
H^3(K,\Q/\Z(2))\to\bigoplus_{v\in X^{(1)}} H^3(K_v,\Q/\Z(2))\stackrel\Sigma\longrightarrow \Q / \Z
$$
is a complex. This also follows from the generalized Weil reciprocity law (\cite{cogal}, Chapter II, Annexe, \S 3) which is valid over an arbitrary field. Therefore the pairing (\ref{nrpairing}) annihilates the diagonal image of $Y(K)$ by a reciprocity law and also its closure in the product topology by a continuity argument combined with Lemma \ref{lemeval}.
This gives rise to the obstruction to weak approximation defined by elements in $H^3_{\rm nr}(K(Y), \Q/\Z(2))$.

Of course, subgroups of $H^3_{\rm nr}(K(Y), \Q/\Z(2))$ define finer obstructions to weak approximation. The following theorem identifies such a subgroup in the case of tori which in fact constitutes the only obstruction.

\begin{theo}\label{mainrecwa}
Let $T$ be a $K$-torus. There is a homomorphism
$$u: \Sha^2_{\omega}(T') \to H^3_{\rm nr}(K(T),\Q/\Z(2))$$
such that each system $(P_v)$ of local points of $T$ annihilated by $\im(u)$ under the pairing {\rm (\ref{nrpairing})} is in the closure of the diagonal image of $T(K)$ for the product topology.
\end{theo}

\begin{dem} Take a flasque resolution
$$
1\to S\to R\to T\to 1
$$
of $T$ as in (\ref{flasque}). It allows one to consider $R$ as $T$-torsor under $S$, with
cohomology class ${[R] \in H^1(T,S)}$.  Choose a smooth compactification $T^c$  of
$T$.
Since $S$ is flasque, by (\cite{ctflasq}, Theorem~2.2) the $T$-torsor
$R$ extends to a $T^c$-torsor $Y \to T^c$ under $S$
 with associated
cohomology class $[Y] \in H^1(T^c, S)$.
The pairing $S\otimes S'\to \Z(2)[2]$ between $S$ and the dual torus $S'$ defined in Section \ref{one} induces a homomorphism
\begin{equation}\label{blimap}
H^1(K,S') \to H^4(T^c,\Z(2))
\end{equation}
given by $a \mapsto a_{T^c} \cup [Y].$

Next recall (e.g. from \cite{kahn}, Prop. 2.9) that there is a natural map
\begin{equation}\label{kahnmap}
H^4(T^c,\Z(2))\to H^3_{\rm nr}(K(T),\Q/\Z(2))
\end{equation}
defined as follows. The Leray spectral sequence for the change-of-sites map $\alpha:\, T^c_{\rm et}\to T^c_{\rm Zar}$ yields an edge map
$$
H^4(T^c, \Z(2))\to H^0_{\rm Zar}(T^c, \RR^4\alpha_*\Z(2)).
$$
According to (\cite{kahn}, Theorem 2.6. c)) the natural map $\Q(2)_{\rm Zar}\to \RR\alpha_*\Q(2)$ is an isomorphism in the derived category of Zariski sheaves.
But since $\Q(2)_{\rm Zar}$ is concentrated in degrees $\leq 2$, we have $\RR^3\alpha_*\Q(2)=\RR^4\alpha_*\Q(2)=0$ and therefore the exact sequence
$$
0\to \Z(i)\to\Q(i)\to \Q/\Z(i)\to 0
$$
induces an isomorphism $\RR^4\alpha_*\Z(2)\cong \RR^3\alpha_*\Q/\Z(2)$. Finally, we have an isomorphism $$ H^0_{\rm Zar}(T^c, \RR^3\alpha_*\Q/\Z(2))\cong H^3_{\rm nr}(K(T), \Q/\Z(2))$$ by the Gersten resolution (see \cite{ctsb}, Theorem 4.1.1).

Observe that although $H^4(T^c, \Z(2))$ and $H^3(T^c, \Q/\Z(2))$ are not isomorphic in general, there is an isomorphism between  $H^4(K(T^c), \Z(2))$ and $H^3(K(T^c), \Q/\Z(2))$ by a similar argument as above. Moreover, the construction of the Gersten resolution implies that the restriction map $H^4(T^c, \Z(2))\to H^4(K(T),\Z(2))$ and the map (\ref{kahnmap}) fit in a commutative diagram
\begin{equation}\label{nrdiag}
\begin{CD}
H^4(T^c, \Z(2)) @>>> H^3_{\rm nr}(K(T), \Q/\Z(2))\\
@VVV @VVV \\
H^4(K(T), \Z(2)) @>\cong>> H^3(K(T), \Q/\Z(2)).
\end{CD}
\end{equation}

Now consider the dual
exact sequence
$$ 1 \to T' \to R' \to S' \to 1$$
coming from the flasque resolution. It induces the exact sequence
$$0 \to H^1(K,S') \to H^2(K,T') \to H^2(K,R'), $$
and since $\Sha^2_{\omega}(R')=0$ by Lemma~\ref{finishaS} (the dual of a quasi-trivial torus being quasi-trivial),
we obtain an isomorphism $\Sha^1_{\omega}(S') \stackrel\sim\to
\Sha^2_{\omega}(T')$.

All in all, following the composite map
$$
\Sha^2_{\omega}(T')\stackrel\sim\to \Sha^1_{\omega}(S')\hookrightarrow H^1(K, S')
$$
by the composition of the maps (\ref{blimap}) and (\ref{kahnmap}) we obtain a homomorphism $u : \Sha^2_{\omega}(T') \to H^3_{\rm nr}(K(T),\Q/\Z(2))$ as in the statement of the theorem.

Let $(P_v) \in \prod_{v \in X^{(1)}} T(K_v)$. Observe that the
evaluation $$ P_v \mapsto [Y](P_v)=[R](P_v) \in H^1(K_v,S)$$ is given by the
coboundary map $\partial_v:\,T(K_v) \to H^1(K_v,S)$ corresponding to the
flasque resolution (\ref{flasque}).
Consider the commutative
diagram with exact rows
$$
\begin{CD}
R(K) @>>> T(K) @>>> H^1(K, S) @>>> 0 \\
@VVV @VVV @VVV \\
\prod_{v \in X^{(1)}} R(K_v) @>>>\prod_{v \in X^{(1)}} T(K_v) @>{(\partial_v)}>> \prod_{v \in X^{(1)}}
H^1(K_v,S) @>>> 0 \cr
&& @VVV @VVV \cr
&& \Sha^2_{\omega}(T')^D @>\simeq>> \Sha^1_{\omega}(S')^D
\end{CD}
$$
in which the vertical maps are induced by the local duality pairings and the horizontal ones by the long exact sequence coming from the flasque resolution. Since $R$ is quasi-trivial, the upper left vertical map has dense image. Hence the diagram together with Theorem~\ref{mainwa} $b)$ shows that $(P_v)$ is in the closure
of $T(K)$ if and only if $(\partial_v(P_v))$ is orthogonal to $\Sha^1_{\omega}(S')$, which means
$$
0=\sum_{v \in X^{(1)}} \langle a_v, \partial_v(P_v) \rangle_v=
\sum_{v \in X^{(1)}}a_v\cup \partial_v(P_v) =\sum_{v \in X^{(1)}} (a_{T^c} \cup [Y])(P_v))
$$
for every $a \in \Sha^1_{\omega}(S')$. Here $(a_{T^c} \cup [Y])(P_v)$ denotes the image of the class ${(a_{T^c} \cup [Y])}\in H^4(T^c,\Z(2))$ by  the evaluation map $H^4(T^c,\Z(2))\to H^4(K_v,\Z(2))\cong\Q/\Z$ associated with  $P_v$; it is 0 for all but finitely many $v$ as $a \in \Sha^1_{\omega}(S')$. Finally, diagram (\ref{nrdiag}) shows that the condition $$\sum_{v \in X^{(1)}} (a_{T^c} \cup [Y])(P_v))=0 $$
means precisely that $(P_v)$ is orthogonal to the image
of $u$ for the pairing (\ref{nrpairing}).
\end{dem}

\begin{remas}\rm ${}$\smallskip

\noindent 1. In their recent work \cite{blimerk}, Blinstein and Merkurjev also construct
a homomorphism
\begin{equation}\label{blimap2}
H^1(K,S') \to H^3_{\nr}(K(T),\Q/\Z(2))
\end{equation}
as a consequence of their Theorems 4.4 and 5.6. By comparing the two constructions (and taking in particular the proof of Proposition 5.3 of \cite{blimerk} into account), one sees that the map of Blinstein and Merkurjev is none but the composite of our maps (\ref{blimap}) and (\ref{kahnmap}).\smallskip

\noindent 2. One can also construct the map (\ref{blimap2}) at a finite level, avoiding the use of the complex $\Z(2)$. Namely, we can find an integer $n>0$ annihilating the group $H^1(K, S')$ (for example, the degree of an extension splitting $S$). We can then send $[Y]\in H^1(T^c, S)$ to a class in $H^2(T^c, {}_nS)$ by a coboundary map $\partial_n$ and lift $a\in H^1(K, S')$ to $a_n\in H^1(K, {}_nS')$. The cup-product $(a_n)_{T^c}\cup\partial_n([Y])$ induced by the pairing ${}_nS'\otimes {}_nS\to \mu_n^{\otimes 2}$ lives in $H^3(T^c, \mu_n^{\otimes 2})$ and restricts to a class in $H^3_{\rm nr}(K(T), \mu_n^{\otimes 2})$ via an argument involving the Gersten resolution as in the construction of the map (\ref{kahnmap}). One can check that this construction is independent of the choices made and yields the same map $H^1(K,S') \to H^3_{\nr}(K(T),\Q/\Z(2))$ as above.\smallskip

\noindent 3. In the classical case where $T$ is defined not over $K$ but over a number field  $k$, we have an isomorphism of cohomology groups $$H^2_{\rm nr}(k(T), \Q/\Z(1))\cong H^2(T^c, \Q/\Z(1))$$ (see e.g. \cite{ctsb}, Theorem 4.2.3). Therefore the evaluation pairing can be defined directly by evaluating classes in $H^2(T^c, \Q/\Z(1))$, without passing through a local ring. In the above proof we have used an intermediate evaluation pairing with the group $H^4(T^c, \Z(2))$ and then composed with the map (\ref{kahnmap}) which is in fact surjective but not injective in general (\cite{kahn}, Proposition 2.9).\smallskip

\noindent 4. Using the same tools as in Section 6 of \cite{pchp}, one can show that the reciprocity obstruction to weak approximation associated with $H^3_{\rm nr}(K(G),\Q/\Z(2))$ is still the only one for a quasi-split reductive $K$-group $G$ such that the universal cover $G^{\rm sc}$ of its derived subgroup satisfies weak approximation. This is in particular the case when $G^{\rm sc}$ is $K$-rational, e.g. $G^{\rm sc}={\rm SL}_n$. The proof requires a non-commutative generalization of Theorem \ref{mainwa} which can be obtained by adapting the methods of \cite{sansuc} in their full generality. Also, one has to use the noncommutative flasque resolutions of \cite{ctresflasq} to prove a non-commutative generalization of Theorem \ref{mainrecwa}.
\end{remas}

\end{document}